\documentclass[11pt,a4paper]{article}
\usepackage{graphicx, amssymb, amsmath, fullpage, times}
\usepackage[bookmarks=false,plainpages=true,pdfstartview=FitBH,
hyperfigures=true,colorlinks=true,linkcolor=blue]{hyperref}
\usepackage{mathrsfs,bm}
\usepackage{slashbox, pifont}

\newtheorem{theorem}{Theorem}
\def\dd#1{\,\mbox{d}#1}

\newtheorem{corollary}[theorem]{Corollary}

\newtheorem{definition}{Definition}

\makeatletter
\def\@yproof[#1]{\@proof{ #1}}
\def\@proof#1{\begin{trivlist}\item[]{\em Proof#1.}}
\newenvironment{proof}{\@ifnextchar[{\@yproof}{\@proof{}}}{
\hfill$\Box$\end{trivlist}\makeatother}

\title{The connectivity-profile of random increasing $k$-trees}

\author{Alexis Darrasse\footnotemark[1]\\
        APR - LIP6\\UPMC\\75005 Paris France\\alexis.darrasse@lip6.fr \and
        Hsien-Kuei Hwang\\
        Institute of Statistical Science\\Academia Sinica\\
        Taipei 115 Taiwan\\hkhwang@stat.sinica.edu.tw \and
        Olivier Bodini\footnotemark[1]\\
        APR - LIP6\\UPMC\\75005 Paris France\\olivier.bodini@lip6.fr \and
        Mich\`ele Soria\thanks{This work was partially supported
        by ANR under the
        contract GAMMA, n\textordmasculine BLAN07-2\_195422.}\\
        APR - LIP6\\UPMC\\75005 Paris France\\michele.soria@lip6.fr}

\date{\today}

\begin{document}

\maketitle

\thispagestyle{empty}

\begin{abstract}
Random increasing $k$-trees represent an interesting, useful class
of strongly dependent graphs for which analytic-combinatorial tools
can be successfully applied. We study in this paper a notion called
connectivity-profile and derive asymptotic estimates for it; some
interesting consequences will also be given.
\end{abstract}


\section{Introduction}
A $k$-tree is a graph reducible to a $k$-clique by successive
removals of a vertex of degree $k$ whose neighbors form a
$k$-clique. This class of $k$-trees has been widely studied in
combinatorics (for enumeration and characteristic
properties~\cite{beineke_number_1969,rose_simple_1974}), in graph
algorithms (many NP-complete problems on graphs can be solved in
polynomial time on $k$-trees~\cite{arnborg_complexity_1987}), and in
many other fields where $k$-trees were naturally encountered
(see~\cite{arnborg_complexity_1987}). By construction, vertices in
such structures are remarkably close, reflecting a highly strong
dependent graph structure, and they exhibit with no surprise the
scale-free property~\cite{gao_degree_2009}, yet somewhat
unexpectedly many properties of random $k$-trees can be dealt with
by standard combinatorial, asymptotic and probabilistic tools, thus
providing an important model of synergistic balance between
mathematical tractability and the predictive power for
practical-world complex networks.

While the term ``$k$-trees" is not very informative and may indeed
be misleading to some extent, they stand out by their underlying
tree structure, related to their recursive definition, which
facilitates the analysis of the properties and the exploration of
the structure.  Indeed, for $k=1$, $k$-trees are just trees, and for
$k\ge 2$ a bijection~\cite{darrasse_limiting_2009} can be explicitly
defined between $k$-trees and a non trivial simple family of trees.

\begin{figure}[!h]
\begin{center}
\begin{tabular}{|c|c|c|c|c|}
\hline
\ding{172} & \ding{173} & \ding{174} & \ding{175} & \ding{176} \\
\includegraphics{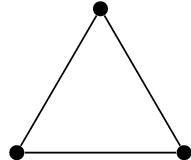} &
\includegraphics{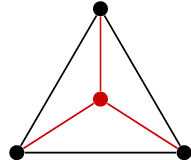} &
\includegraphics{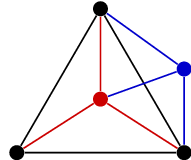} &
\includegraphics{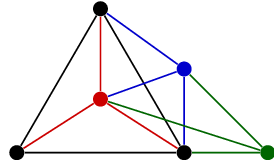} &
\includegraphics{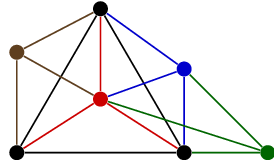} \\
\includegraphics{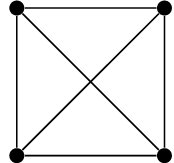} &
\includegraphics{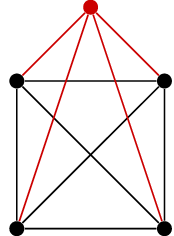} &
\includegraphics{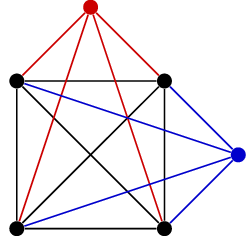} &
\includegraphics{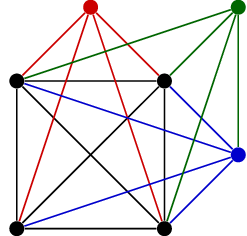} &
\includegraphics{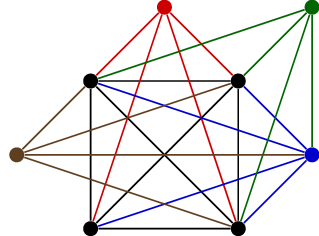} \\ \hline
\end{tabular}
\end{center}
\caption{\emph{The first few steps of generating a $3$-tree and a
$4$-tree. Obviously, these graphs show the high connectivity of
$k$-trees.}} \label{fg-2-tr}
\end{figure}

The process of generating a $k$-tree begins with a $k$-clique, which
is itself a $k$-tree; then the $k$-tree grows by linking a new
vertex to every vertex of an existing $k$-clique, and to these
vertices only. The same process continues; see Figure~\ref{fg-2-tr}
for an illustration. Such a simple process is reminiscent of several
other models proposed in the literature such as
$k$-DAGs~\cite{devroye_long_2009},  random
circuits~\cite{arya_expected_1999}, preferential
attachment~\cite{barabsi_emergence_1999,bollobs_degree_2001,
hwang_profiles_2007}, and many other models (see, for
example,~\cite{boccaletti_complex_2006,
durrett_random_2006,newman_structure_2003}). While the construction
rule in each of these models is very similar, namely, linking a new
vertex to $k$ existing ones, the mechanism of choosing the existing
$k$ vertices differs from one case to another, resulting in very
different topology and dynamics.

Restricting to the procedure of choosing a $k$-clique each time a
new vertex is added, there are several variants of $k$-trees
proposed in the literature depending on the modeling needs. So
$k$-trees can be either labeled~\cite{beineke_number_1969},
unlabeled~\cite{labelle_labelled_2004},
increasing~\cite{zhang_high-dimensional_2006},
planar~\cite{zhang_high-dimensional_2006},
non-planar~\cite{beineke_number_1969}, or
plane~\cite{palmer_number_1973}, etc.

For example, the family of random Apollonian networks, corresponding
to planar 3-trees, has recently been employed as a model for complex
networks~\cite{andrade_apollonian_2005,zhang_high-dimensional_2006}.
In these frameworks, since the exact topology of the real networks
is difficult or even impossible to describe, one is often led to the
study of models that present similarities to some observed
properties such as the degree of a node and the distance between two
nodes of the real structures.

For the purpose of this paper, we distinguish between two models of
random labeled non-plane $k$-trees; by non-plane we mean that we
consider these graphs as given by a set of edges (and not by its
graphical representation):
\begin{itemize}

\item[--] \emph{random simply-generated $k$-trees}, which correspond
to a uniform probability distribution on this class of $k$-trees,
and

\item[--] \emph{random increasing $k$-trees}, where we consider the
iterative generation process:  at each time step, all existing
$k$-cliques are equally likely to be selected and the new vertex is
added with a label which is greater than the existing ones.

\end{itemize}
The two models are in good analogy to the simply-generated family of
trees of Meir and Moon~\cite{meir_altitude_1978} marked specially by
the functional equation $f(z) = z\Phi(f(z))$ for the underlying
enumerating generating function, and the increasing family of trees
of Bergeron et al.~\cite{bergeron_varieties_1992}, characterized by
the differential equation $f'(z) = \Phi(f(z))$. Very different
stochastic behaviors have been observed for these families of trees.
While similar in structure to these trees, the analytic problems on
random $k$-trees we are dealing with here are however more involved
because instead of a scalar equation (either functional, algebraic,
or differential), we now have a system of equations.

\begin{table}
\begin{center}
\begin{tabular}{|r||c|c|} \hline
\backslashbox{Properties}{Model} & Simply-generated structures&
Increasing structures \\ \hline Combinatorial description &
$\mathcal{T}_s = \mbox{Set}(\mathcal{Z}\times\mathcal{T}_s^k)$ &
$\mathcal{T} =
\mbox{Set}(\mathcal{Z}^\square \times\mathcal{T}^k)$\\
\hline Generating function & $T_s(z) = \exp(z T_s^k(z))$
& $ T'(z) =  T^k(z)$ \\
\hline Expansion near singularity& $T_s(z) = \tau -
h\sqrt{1-z/\rho}+\ldots$ & $ T(z) = (1-kz)^{-1/k}$ \\ \hline Mean
distance of nodes & $O(\sqrt{n})$ & $O(\log n)$ \\ \hline Degree
distribution & Power law with  exp.\ tails & Power
law~\cite{gao_degree_2009}
\\\hline Root-degree distribution & Power law with  exp.\ tails &
Stable law~(Theorem~\ref{thm-ld}) \\\hline Expected Profile &
Rayleigh limit law & Gaussian limit law (\ref{Ecp-LLT})
\\ \hline
\end{tabular}
\end{center}
\caption{\emph{The contrast of some properties between
random simply-generated $k$-trees and
random increasing $k$-trees. Here $\mathcal{Z}$ denotes a node
and $\mathcal{Z}^\square$ means a marked node.}} \label{tb1}
\end{table}

It is known that random trees in the family of increasing trees are
often less skewed, less slanted in shape, a typical description
being the logarithmic order for the distance of two randomly chosen
nodes; this is in sharp contrast to the square-root order for random
trees belonging to the simply-generated family; see for
example~\cite{bergeron_varieties_1992,drmota_random_2009,
fuchs_profiles_2006,marckert_families_2008,meir_altitude_1978}. Such
a contrast has inspired and stimulated much recent research. Indeed,
the majority of random trees in the literature of discrete
probability, analysis of algorithms, and random combinatorial
structures are either $\log n$-trees or $\sqrt{n}$-trees, $n$ being
the tree size. While the class of $\sqrt{n}$-trees have been
extensively investigated by probabilists and combinatorialists,
$\log n$-trees are comparatively less addressed, partly because most
of them were encountered not in probability or in combinatorics, but
in the analysis of algorithms.

Table~\ref{tb1} presents a comparison of the two models: the classes
${\mathcal{T}}_s$ and $\mathcal{T}$, corresponding respectively to
simply-generated $k$-trees and increasing $k$-trees. The results
concerning simple $k$-trees are given
in~\cite{darrasse_limiting_2009, darrasse_unifying_????}, and those
concerning increasing $k$-trees are derived in this paper (except
for the power law distribution~\cite{gao_degree_2009}). We start
with the specification, described in terms of operators of the
symbolic method~\cite{flajolet_analytic_2009}. A structure of
${\mathcal{T}}_s$ is a set of $k$ structures of the same type, whose
roots are attached to a new node: $\mathcal{T}_s =
\mbox{Set}(\mathcal{Z}\times\mathcal{T}_s^k)$, while a structure of
${\mathcal{T}}$ is an increasing structure, in the sense that the
new nodes get labels that are smaller than those of the underlying
structure (this constraint is reflected by the box-operator)
$\mathcal{T} = \mbox{Set} (\mathcal{Z}^\square\times\mathcal{T}^k)$.
The analytic difference immediately appears in the enumerative
generating functions that translate the specifications: the
simply-generated structures are defined by $T_s(z)= \exp(z
T_s^k(z))$ and corresponding increasing structures satisfy the
differential equation $ T'(z)= T^k(z)$. These equations lead to a
singular expansion of the square-root type in the simply-generated
model, and a singularity in $(1 - kz)^{-1/k}$ in the increasing
model. Similar analytic differences arise in the bivariate
generating functions of shape parameters.

The expected distance between two randomly chosen vertices or the
average path length is one of the most important shape measures in
modeling complex networks as it indicates roughly how efficient the
information can be transmitted through the network. Following the
same $\sqrt{n}$-vs-$\log n$ pattern, it is of order $\sqrt n$ in the
simply-generated model, but $\log n$ in the increasing model.
Another equally important parameter is the degree distribution of a
random vertex: its limiting distribution is a power law with
exponential tails in the simply-generated model of the form
$d^{-3/2} \rho_k^d$, in contrast to a power-law in the increasing
model of the form $d^{-1-k/(k-1)}$, $d$ denoting the
degree~\cite{gao_degree_2009}. As regards the degree of the root,
its asymptotic distribution remains the same as that of any vertex
in the simply-generated model, but in the increasing model, the
root-degree distribution is different, with an asymptotic stable law
(which is Rayleigh in the case $k=2$); see Theorem~\ref{thm-ld}.

Our main concern in this paper is the connectivity-profile.
Recall that the profile of an usual tree is the sequence
of numbers, each enumerating the total number of nodes with the same
distance to the root. For example, the tree
\includegraphics{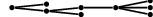}
has the profile $\{1,2,2,1,3\}$. Profiles represent one of the
richest shape measures and they convey much information regarding
particularly the silhouette. On random trees, they have been
extensively studied recently; see~\cite{chauvin_martingales_2005,
drmota_profile_1997,drmota_functional_2008,fuchs_profiles_2006,
hwang_profiles_2007,marckert_families_2008,park_profiles_2009}.
Since $k$-trees have many cycles for $k\ge2$, we call the profile of
the transformed tree (see next section) \emph{the
connectivity-profile} as it measures to some extent the connectivity
of the graph. Indeed this connectivity-profile corresponds to the
profile of the ``shortest-path tree'' of a $k$-tree, as defined by
Proskurowski~\cite{proskurowski_k-trees:_1980}, which is nothing
more than the result of a Breadth First Search (BFS) on the graph.
Moreover, in the domain of complex networks, this kind of BFS trees
is an important object; for example, it describes the results of the
\texttt{traceroute} measuring
tool~\cite{stevens_chapter_1994,viger_detection_2008} in the study
of the topology of the Internet.

We will derive precise asymptotic approximations to the expected
connectivity-profile of random increasing $k$-trees, the major tools
used being based on the resolution of a system of differential
equations of Cauchy-Euler type (see~\cite{chern_asymptotic_2002}).
In particular, the expected number of nodes at distance $d$ from the
root follows asymptotically a Gaussian distribution, in contrast to
the Rayleigh limit distribution in the case of simply-generated
$k$-trees. Also the limit distribution of the number of nodes with
distance $d$ to the root will be derived when $d$ is bounded. Note
that when $d=1$, the number of nodes at distance $1$ to the root is
nothing but the degree of the root.

This paper is organized as follows. We first present the definition
and combinatorial specification of random increasing $k$-trees in
Section~\ref{sec:def}, together with the enumerative generating
functions, on which our analytic tools will be based. We then
present two asymptotic approximations to the expected
connectivity-profile in Section~\ref{sec:profile}, one for $d=o(\log
n)$ and the other for $d\to\infty$ and $d=O(\log n)$. Interesting
consequences of our results will also be given. The limit
distribution of the connectivity-profile in the range when $d=O(1)$
is then given in Section~\ref{sec:ld}.

\section{Random increasing $k$-trees and generating functions}
\label{sec:def}

Since $k$-trees are graphs full of cycles and cliques, the key step
in our analytic-combinatorial approach is to introduce a bijection
between $k$-trees and a suitably defined class of trees (\emph{bona
fide} trees!) for which generating functions can be derived. This
approach was successfully applied to simply-generated family of
$k$-trees in~\cite{darrasse_limiting_2009}, which leads to a system
of algebraic equations. The bijection argument used there can be
adapted \emph{mutatis mutandis} here for increasing $k$-trees, which
then yields a system of differential equations through the bijection
with a class of increasing trees~\cite{bergeron_varieties_1992}.

\begin{figure}
\begin{center}
\includegraphics{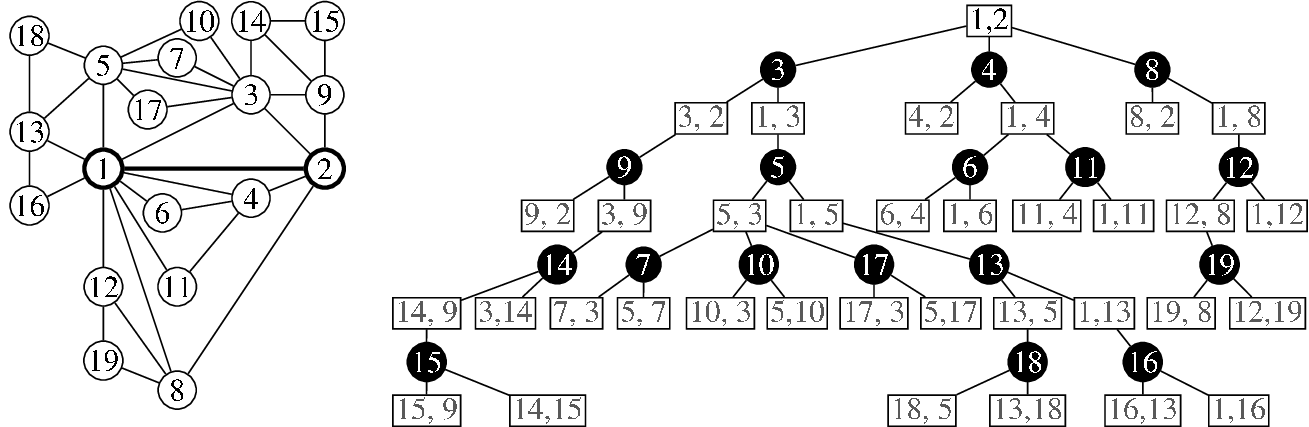}
\end{center}
\caption{\emph{A $2$-tree (left) and its corresponding increasing
tree representation (right).}}\label{fig:bij}
\end{figure}

\paragraph{Increasing $k$-trees and the bijection.}
Recall that a $k$-clique is a set of $k$ mutually adjacent vertices.
\begin{definition}
An increasing $k$-tree is defined recursively as follows. A
$k$-clique in which each vertex gets a distinct label from
$\{1,\dots,k\}$ is an increasing $k$-tree of $k$ vertices. An
increasing $k$-tree with $n > k$ vertices is constructed from an
increasing $k$-tree with $n-1$ vertices by adding a vertex labeled
$n$ and by connecting it by an edge to each of the $k$ vertices in
an existing $k$-clique.
\end{definition}

By \emph{random increasing $k$ trees}, we assume that all existing
$k$-cliques are equally likely each time a new vertex is being
added. One sees immediately that the number $T_n$ of increasing
$k$-trees of $n+k$ nodes is given by $T_n = \prod_{0\le i<n}(ik+1)$.

Note that if we allow any permutation on all labels, we obtain the
class of simply-generated $k$-trees where monotonicity of labels
along paths fails in general.

Combinatorially, simply-generated $k$-trees are in
bijection~\cite{darrasse_limiting_2009} with the family of trees
specified by $\mathcal{K}_s = \mathcal{Z}^k \times \mathcal{T}_s$,
where $\mathcal{T}_s = \mbox{Set}(\mathcal{Z} \times
\mathcal{T}_s^k)$. Given a rooted $k$-tree $G$ of $n$ vertices, we
can transform $G$ into a tree $T$, with the root node labeled
$\{1,\dots,k\}$, by the following procedure. First, associate a
white node to each $k$-clique of $G$ and a black node to each
$(k+1)$-clique of $G$. Then add a link between each black node and
all white nodes associated to the $k$-cliques it contains. Each
black node is labeled with the only vertex not appearing in one of
the black nodes above it or in the root.  The last step in order to
complete the bijection is to order the $k$ vertices of the root and
propagate this order to the $k$ sons of each black node. This
constructs a tree from a $k$-tree (see Figure~\ref{fig:bij});
conversely, we can obtain the $k$-tree through a simple traversal of
the tree.

Such a bijection translates directly to increasing $k$-trees by
restricting the class of corresponding trees to those respecting a
monotonicity constraint on the labels, namely, on any path from the
root to a leaf the labels are in increasing order. This yields the
combinatorial specification of the class of increasing trees
$\mathcal{T} = \mbox{Set}(\mathcal{Z}^\square \times
\mathcal{T}^k)$. An increasing $k$-tree is just a tree in
$\mathcal{T}$ together with the sequence $\{1,\dots,k\}$
corresponding to the labels of the root-clique\footnote{We call
\textit{root-clique} the clique composed by the $k$ vertices
$(1,\ldots,k)$. The increasing nature of the $k$-trees guarantees
that these vertices always form a clique. We call
\textit{root-vertex} the vertex with label $1$.}. A tree in
$\mathcal{K}$ is thus completely determined by its $\mathcal{T}$
component, giving $\mathcal{K}_{n+k} \equiv \mathcal{T}_n$. For
example figure~\ref{fig:bij} shows a $2$-tree with $19$ vertices and
its tree representation with $17$ black nodes. In the rest of this
paper we will thus focus on class $\mathcal{T}$.

\paragraph{Generating functions.}
Following the bijection, we see that the complicated dependence
structure of $k$-trees is now completely described by the class of
increasing trees specified by $\mathcal{T} =
\mbox{Set}(\mathcal{Z}^\square \times \mathcal{T}^k)$. For example,
let $T(z) := \sum_{n\ge0} T_n z^n/n!$ denote the exponential
generating function of the number $T_n$ of increasing $k$-trees of
$n+k$ vertices. Then the specification translates into the equation
\[
    T(z) = \exp\left(\int_0^z T^k(x) \dd x\right),
\]
or, equivalently, $T'(z) = T^{k+1}(z)$ with $T(0)=1$, which is
solved to be
\[
    T(z) = (1-kz)^{-1/k},
\]
we then check that $T_n = \prod_{0\le i<n}(ik+1)$.

If we mark the number of neighbors of the root-node in $\mathcal{T}$
by $u$, we obtain
\[
    T(z,u) = \exp\left(u \int_0^z T(x) T^{k-1}(x,u) \dd x\right),
\]
where the coefficients $n![u^\ell z^n] T(z,u)$ denote the number of
increasing $k$-trees of size $n+k$ with root degree equal to
$k+\ell-1$. Taking derivative with respect to $z$ on both sides and
then solving the equation, we get the closed-form expression
\begin{align} \label{F-sol}
    T(z,u) = \left(1-u(1-(1-kz)^{1-1/k})\right)^{-1/(k-1)}.
\end{align}

Since $k$-trees can be transformed into ordinary increasing trees,
the profiles of the transformed trees can be naturally defined,
although they do not correspond to simple parameters on $k$-trees.
While the study of profiles may then seem artificial, the results do
provide more insight on the structure of random $k$-trees. Roughly,
we expect that all vertices on $k$-trees are close, one at most of
logarithmic order away from the other. The fine results we derive
provide in particular an upper bound for that.

Let $X_{n;d,j}$ denote the number of nodes at distance $d$ from $j$
vertices of the root-clique in a random $k$-tree of $n+k$ vertices.
Let $T_{d,j}(z,u)= \sum_{n\ge0} T_n \mathbb{E}(u^{X_{n;d,j}})
z^n/n!$ denote the corresponding bivariate generating function.
\begin{theorem} The generating functions $T_{d,j}$'s satisfy the
differential equations
\begin{equation} \label{bgf-Tdj}
    \frac{\partial}{\partial z} T_{d,j}(z,u) =
    u^{\delta_{d,1}}T_{d,j-1}^j(z,u) T_{d,j}^{k-j+1}(z,u),
\end{equation}
with the initial conditions $T_{d,j}(0,u)=1$ for $1\le j\le k$,
where $\delta_{a,b}$ denotes the Kronecker function,
$T_{0,k}(z,u)=T(z)$ and $T_{d,0}(z,u) = T_{d-1,k}(z,u)$.
\end{theorem}
\begin{proof}\ The theorem follows from
\begin{equation*}
    T_{d,j}(z,u) = \exp\left(u^{\delta_{d,1}}
    \int_0^z T_{d,j-1}^j(x,u) T_{d,j}^{k-j}(x,u)\dd x\right),
\end{equation*}
with $T_{d,j}(z,1) = T(z)$.
\end{proof}

For operational convenience, we normalize all $z$ by $z/k$ and write
$\tilde{T}(z) := T(z/k) = (1-z)^{-1/k}$. Similarly, we define
$\tilde{T}_{d,j}(z,u) := T_{d,j}(z/k,u)$ and have, by
(\ref{bgf-Tdj}),
\begin{align} \label{Tdj}
    \frac{\partial}{\partial z} \tilde{T}_{d,j}(z,u)
    =\frac{u^{\delta_{d,1}}}{k}\tilde{T}_{d,j-1}^j(z,u)
    \tilde{T}_{d,j}^{k-j+1}(z,u),
\end{align}
with $\tilde{T}_{d,j}(1,z) = \tilde{T}(z)$,
$\tilde{T}_{0,k}(z,u)=\tilde{T}(z)$ and $\tilde{T}_{d,0}(z,u) =
\tilde{T}_{d-1,k}(z,u)$.

\section{Expected connectivity-profile}
\label{sec:profile}

We consider the expected
connectivity-profile $\mathbb{E}(X_{n;d,j})$ in this section.
Observe first that
\[
    \mathbb{E}(X_{n;d,j}) = \frac{k^n[z^n]\tilde{M}_{d,j}(z)}
    {T_n} ,
\]
where $\tilde{M}_{d,j}(z) := \partial \tilde{T}_{d,j}(z,u)/(\partial
u)|_{u=1}$. It follows from (\ref{Tdj}) that
\begin{align}\label{Mdj}
    \tilde{M}_{d,j}'(z) = \frac1{k(1-z)}
    \left((k-j+1)\tilde{M}_{d,j}(z)
    +j\tilde{M}_{d,j-1}(z) + \delta_{d,1}\tilde{T}(z)\right).
\end{align}
This is a standard differential equation of Cauchy-Euler type whose
solution is given by (see~\cite{chern_asymptotic_2002})
\[
    \tilde{M}_{d,j}(z) = \frac{(1-z)^{-(k-j+1)/k}}k
    \int_0^z (1-x)^{-(j-1)/k}\left( j\tilde{M}_{d,j-1}(x)
    +\delta_{d,1}\tilde{T}(x) \right) \dd x,
\]
since $\tilde{M}_{d,j}(0)=0$. Then, starting from
$\tilde{M}_{0,k}=0$, we get
\begin{align*}
    \tilde{M}_{1,1}(z)= \frac1{k-1}\left(\frac1{1-z}
    - \frac1{(1-z)^{1/k}}\right)
    = \frac{\tilde{T}^k(z) - \tilde{T}(z)}{k-1}.
\end{align*}
Then by induction, we get
\[
    \tilde{M}_{d,j}(z)
    \sim \frac{j}{(k-1)(d-1)!}\cdot\frac{1}{1-z}
    \log^{d-1}\frac1{1-z} \qquad(1\le j\le k;d\ge1; z\sim1).
\]
So we expect, by singularity analysis, that
\[
    \mathbb{E}(X_{n;d,j}) \sim \Gamma(1/k)\frac{j}{k-1}
    \cdot \frac{(\log n)^{d-1}}{(d-1)!}\,n^{1-1/k},
\]
for large $n$ and fixed $d$, $k$ and $1\le j\le k$. We can indeed
prove that the same asymptotic estimate holds in a larger range.

\begin{theorem} \label{thm-E}
The expected connectivity-profile $\mathbb{E}(X_{n;d,j})$
satisfies for $1\le d= o(\log n)$
\begin{align}\label{Ecp-1}
     \mathbb{E}(X_{n;d,j}) \sim \Gamma(1/k)\frac{j}{k-1}
     \cdot \frac{(\log n)^{d-1}}{(d-1)!}\,n^{1-1/k},
\end{align}
uniformly in $d$, and for $d\to\infty$, $d=O(\log n)$,
\begin{align}\label{Ecp2}
    \mathbb{E}(X_{n;d,j}) \sim
    \frac{\Gamma(1/k)h_{j,1}(\rho)
    \rho^{-d} n^{\lambda_1(\rho)-1/k}}
    {\Gamma(\lambda_1(\rho))\sqrt{2\pi(\rho\lambda_1'(\rho)
    +\rho^2\lambda_1''(\rho))\log n}}
\end{align}
where $\rho=\rho_{n,d}>0$ solves the equation $\rho\lambda_1'(\rho)=
d/\log n$, $\lambda_1(w)$ being the largest zero (in real part) of
the equation $\prod_{1\le \ell\le k}(\theta-\ell/k)- k! w/k^k=0$ and
satisfies $\lambda_1(1) =(k+1)/k$.
\end{theorem}
An explicit expression for the $h_{j,1}$'s is given as follows. Let
$\lambda_1(w),\ldots,\lambda_k(w)$ denote the zeros of the equation
$\prod_{1\le \ell\le k}(\theta-\ell/k)- k! w/k^k=0$. Then for $1\le
j\le k$
\begin{align}\label{hj1}
    h_{j,1}(w)  = \frac{j!w(w-1)}{(k\lambda_1(w)-1)\left(
    \sum_{1\le s\le k}\frac1{k\lambda_1(w)-s}\right)
    \prod_{k-j+1\le s\le k+1}(k\lambda_1(w)-s)}.
\end{align}

The theorem cannot be proved by the above inductive argument and our
method of proof consists of the following steps. First, the
bivariate generating functions $\mathscr{M}_j(z,w) := \sum_{d\ge1}
\tilde{M}_{d,j}(z) w^d$ satisfy the linear system
\[
    \left((1-z)\frac{\mbox{d}}{\dd{z}}
    - \frac{k-j+1}{k}\right)\mathscr{M}_j
    = \frac{j}{k} \mathscr{M}_{j-1}
    + \frac{w\tilde{T}}{k}\qquad(1\le j\le k).
\]
Second, this system is solved and has the solutions
\[
    \mathscr{M}_j(z,w) = \sum_{1\le j\le k}
    h_{j,m}(w)(1-z)^{-\lambda_m(w)}
    - \frac{w-(w-1)\delta_{k,j}}{k}\,\tilde{T}(z),
\]
where the $h_{j,m}$ have the same expression as $h_{j,1}$ but with
all $\lambda_1(w)$ in (\ref{hj1}) replaced by $\lambda_m(w)$. While
the form of the solution is well anticipated, the hard part is the
calculations of the coefficient-functions $h_{j,m}$. Third, by
singularity analysis and a delicate study of the zeros, we then
conclude, by saddle-point method, the estimates given in the
theorem.

\begin{corollary} The expected degree of the root
$\mathbb{E}(X_{n,1,j})$ satisfies
\[
    \mathbb{E}(X_{n,1,j}) \sim \Gamma(1/k)\frac{j}{k-1}
     \,n^{1-1/k}\qquad(1\le j\le k).
\]
\end{corollary}
This estimate also follows easily from (\ref{F-sol}).

Let $H_k := \sum_{1\le\ell\le k} 1/\ell$ denote the harmonic
numbers and $H_k^{(2)} := \sum_{1\le\ell\le k} 1/\ell^2$.
\begin{corollary} The expected number of nodes at distance
$d= \left\lfloor \frac{1}{kH_k}\log n + x\sigma\sqrt{\log
n}\right\rfloor$ from the root, where $\sigma =
\sqrt{H_k^{(2)}/(kH_k^3)}$, satisfies, uniformly for $x=o((\log
n)^{1/6})$,
\begin{align}\label{Ecp-LLT}
    \mathbb{E}(X_{n;d,j})\sim
    \frac{n e^{-x^2/2}}{\sqrt{2\pi\sigma^2\log n}}.
\end{align}

\end{corollary}
This Gaussian approximation justifies the last item corresponding to
increasing trees in Table~\ref{tb1}.

Note that $\lambda_1(1)=(k+1)/k$ and $\alpha = d/\log n \sim
1/(kH_k)$. In this case, $\rho=1$ and
\[
    \rho\lambda_1'(\rho) = \frac1{\sum_{1\le \ell\le k}
    \frac{1}{\lambda_1(\rho)-\frac \ell k}},
\]
which implies that $\lambda_1(\rho)-1/k -\alpha\log \rho \sim 1$.

\begin{corollary} \label{cor-height}
Let $\mathscr{H}_{n;d,j} := \max_d X_{n;d,j}$ denote the
height of a random increasing $k$-tree of $n+k$ vertices. Then
\[
    \mathbb{E}(\mathscr{H}_n) \le \alpha_+\log n -
   \frac{\alpha_+}{2(\lambda_1(\alpha_+)-\frac1k)}\log\log n
   + O(1),
\]
where $\alpha_+>0$ is the solution of the system of equations
\[
    \left\{\begin{array}{l}
        \displaystyle\frac{1}{\alpha_+} = \sum_{1\le \ell\le k}
        \frac1{v-\frac \ell k},\\ \displaystyle
        v-\frac1k-\alpha_+\sum_{1\le \ell \le k}
        \log\left(\frac{k}{\ell} v-1\right) = 0.
    \end{array}\right.
\]
\end{corollary}
Table~\ref{tab1} gives the numerical values of $\alpha_+$ for small
values of $k$.
\begin{table}[h!]
\begin{center}
\begin{tabular}{|c||c|c|c|c|c|}\hline
$k$ & $2$ & $3$ & $4$ & $5$ & $6$ \\ \hline%
$\alpha_+$ & $1.085480$ & $0.656285$ & $0.465190$ &
$0.358501$ & $0.290847$\\ \hline\hline%
$k$ & $7$ & $8$ & $9$ & $10$ & $20$ \\ \hline %
$\alpha_+$ & $0.244288$ & $0.210365$ & $0.184587$ & $0.164356$ &
$0.077875$\\ \hline
\end{tabular}
\end{center}
\caption{\emph{Approximate numerical values of $\alpha_+$.}}
\label{tab1}
\end{table}
For large $k$, one can show that $\alpha_+\sim 1/(k\log 2)$ and
$\lambda_1(\alpha_+) \sim 2$.

Corollary~\ref{cor-height} justifies that the mean distance of
random $k$-trees are of logarithmic order in size, as stated in
Table~\ref{tb1}.

\begin{corollary} The width $\mathscr{W}_{n;d,j} := \max_d X_{n;d,j}$
is bounded below by
\[
    \mathbb{E}(\mathscr{W}_n) =
    \mathbb{E}(\max_d X_{n,d}) \ge \max_d \mathbb{E}(X_{n,d})
    \asymp \frac{n}{\sqrt{\log n}}.
\]
\end{corollary}

We may conclude briefly from all these results that \emph{in the
transformed increasing trees of random increasing $k$-trees, almost
all nodes are located in the levels with $d= \frac{1}{kH_k}\log n +
O(\sqrt{\log n})$, each with $n/\sqrt{\log n}$ nodes.}

\section{Limiting distributions}
\label{sec:ld}

With the availability of the bivariate generating functions
(\ref{bgf-Tdj}), we can proceed further and derive the limit
distribution of $X_{n;d,j}$ in the range where $d=O(1)$. The case
when $d\to\infty$ is much more involved; we content ourselves in
this extended abstract with the statement of the result for bounded
$d$.

\begin{theorem} \label{thm-ld}
The random variables $X_{n;d,j}$, when normalized by their mean
orders, converge in distribution to
\begin{align}\label{Xndj-cid}
    \frac{X_{n;d,j}}{n^{1-1/k}(\log n)^{d-1}/(d-1)!}
    \stackrel{d}{\to}\Xi_{d,j},
\end{align}
where
\begin{align*}
    \mathbb{E}(e^{\Xi_{d,j} u}) &=\Gamma(\tfrac1k)
    \sum_{m\ge0} \frac{c_{d,j,m}}{m!\Gamma(m(1-1/k)+1/k)}\,u^m\\
    &= \frac{\Gamma(\frac1k)}{2\pi i}\int_{-\infty}^{(0+)} e^\tau
    \tau^{-1/k}C_{d,j}\left(\tau^{-1+1/k} u\right) \dd \tau,
\end{align*}
and $C_{d,j}(u) :=1+ \sum_{m\ge1} c_{d,j,m} u^m/m!$ satisfies the
system of differential equations
\begin{align} \label{Cdju}
    (k-1)uC_{d,j}'(u) + C_{d,j}(u) = C_{d,j}(u)^{k+1-j}
    C_{d,j-1}(u)^j\qquad(1\le j\le k),
\end{align}
with $C_{d,0}=C_{d-1,k}$. Here the symbol $\int_{-\infty}^{(0+)}$
denotes any Hankel contour starting from $-\infty$ on the real axis,
encircling the origin once counter-clockwise, and returning to
$-\infty$.
\end{theorem}

We indeed prove the convergence of all moments, which is stronger
than weak convergence; also the limit law is uniquely determined by
its moment sequence.

So far only in special cases do we have explicit solution for
$C_{1,j}$: $C_{1,1}(u) = (1+u)^{-1/(k-1)}$ and
\[
    C_{1,2}(u) = \left\{\begin{array}{ll}
    \frac{e^{1/(1+u)}}{1+u},&\text{if } k=2;\\
    \frac1{1+u^{1/2}\arctan(u^{1/2})},&
    \text{if } k=3.
    \end{array}\right.
\]

Note that the result (\ref{Xndj-cid}) when $d=0$ can also be derived
directly by the explicit expression (\ref{F-sol}). In particular,
when $k=2$, the limit law is Rayleigh.

\bibliographystyle{plain}
\bibliography{increasing_ktrees}
\end{document}